\documentclass[12pt]{article}
\usepackage{harvard}
\usepackage{amssymb}
\DeclareMathAlphabet{\mathscr}{OT1}{pzc}{m}{it}
\def\be{\begin{equation}}
\def\ee{\end{equation}}
\def\bea{\begin{eqnarray*}}
\def\eea{\end{eqnarray*}}
\def\RR{\mathbb R}
\def\CP{\mathbb C \mathbb P}

\newtheorem{defn}{Definition}

\newtheorem{thm}{Theorem}
\newtheorem{prop}{Proposition}

\newenvironment{proof}{\medskip \noindent
{\bf Proof.}}{\hfill \rule{.5em}{1em} \medskip}
\begin{document}
\sloppy
\title{The Einstein-Maxwell Equations,\\ Extremal K\"ahler Metrics, \\ and Seiberg-Witten Theory}

\author{Claude LeBrun\thanks{Supported 
in part by  NSF grant DMS-0604735.} 
\\ 
SUNY Stony
 Brook 
  }

\date{}
\maketitle


Nigel Hitchin has played a key role in the exploration
 of $4$-dimensional Riemannian geometry, and in particular 
 (Atiyah, Hitchin \harvardand\ Singer \citeyear*{AHS}, 
 Hitchin \citeyear*{hitharm,hit,hitpos,hitpoly,hitka},
Hitchin, Karlhede, Lindstr{\"o}m \harvardand\
  Ro{\v{c}}ek  \citeyear*{HKLR})
has made foundational  contributions   to the 
theory of self-dual manifolds,  $4$-dimensional 
Einstein manifolds, spin$^c$ structures, and K\"ahler geometry. 
In the process, he has often 
alerted the rest of us to the profound mathematical interest   
 of  beautiful  geometric problems that 
had previously only been considered by physicists. I would therefore like to use
the occasion of 
Nigel's 60$^{\rm th}$ birthday as an opportunity to draw the attention
of an audience of geometers and physicists to 
an interesting relationship between the $4$-dimensional 
Einstein-Maxwell equations and K\"ahler geometry, and 
point out some fascinating open problems that  directly impinge 
on this relationship.

Let us begin by recalling  that 
a $2$-form $F$ on an oriented  pseudo-Riemannian $n$-manifold 
$(M,g)$ is said to satisfy {\em Maxwell's equations} (in vacuo) if and only if
\begin{eqnarray*}
dF&=&0\\
d\star F&=&0
\end{eqnarray*}
where $\star$ is the Hodge star operator. If $M$ is compact and 
$g$ is Riemannian, these equations  of course just mean that 
$F$ is a harmonic $2$-form, and Hodge theory thus asserts  that 
there is in fact exactly one solution in each de Rham cohomology class
$[F]\in H^2 (M, \RR)$. This solution may be found by 
minimizing the action 
$$F \longmapsto  \int_M|F|^2_gd\mu_g$$
among all closed forms $F\in [F]$. In this context, dimension four enjoys a 
somewhat privileged status, because it is precisely when $n=4$  
that both the action  and  the solutions themselves
become conformally invariant, in the sense that they are unaltered by 
replacing $g$ with any conformally related metric $\tilde{g}=u^2g$.

When $n > 2$, coupling  these equations to the gravitational field
\cite{hawkell,mtw} 
gives rise to the so-called {\em Einstein-Maxwell equations}
(with cosmological constant) 
\begin{eqnarray*}
dF&=&0\\
d\star F&=&0\\
\Big[r+ F\circ F\Big]_0&=&0
\end{eqnarray*}
where $r$ is the Ricci tensor of $g$, where 
$(F\circ F)_{jk}= {F_j}^\ell F_{\ell k}$ is obtained by
composing $F$ with itself as an endomorphism of $TM$, and  
where 
$[~]_0$ indicates
the trace-free part with respect to $g$. 
In the compact Riemannian setting, these equations may be understood 
as the Euler-Lagrange equations of the functional 
$$(g,F) \longmapsto   \int_M(s_g+|F|^2_g)d\mu_g$$
where $F$ is again allowed to vary over all
closed $2$-forms in a given de Rham class, and where
$g$ is allowed to vary over all Riemannian metrics 
of some fixed total volume $V$. 
The privileged status of dimension four
becomes  more pronounced in this 
context, for it is only when $n=4$  
that the Einstein-Maxwell equations imply
that the scalar curvature $s$ of $g$ is 
{\em constant}. Indeed, this just reflects 
Yamabe's observation \cite{yamabe} that a Riemannian metric
has constant scalar curvature iff it is a 
critical point of the restriction of the Einstein-Hilbert action 
$\int s~d\mu$ to  the space of  volume-$V$ metrics 
in a fixed   conformal class.
When $n=4$, the
 conformal invariance of $\int |F|^2d\mu$ 
 thus implies that critical points of the above functional 
 must have constant scalar curvature; 
 but when $n\neq 4$, by contrast, 
 the scalar curvature  turns out  to be constant  
 only when    $F$ has constant norm.

We have just observed that
  the Einstein-Maxwell equations on a $4$-manifold imply that  the
scalar curvature is constant. But  what 
happens  in the converse direction is far more remarkable: 
namely, {\em any constant-scalar-curvature K\"ahler metric 
on a $4$-manifold may be interpreted as a solution of the Einstein-Maxwell
equations}. Indeed, suppose that $(M^4,g,J)$ is a K\"ahler surface, with 
K\"ahler form
$\omega = g (J\cdot , \cdot )$
and Ricci form
$\rho= r(J\cdot , \cdot )$.
Let 
$$\mathring{\rho} = \rho - \frac{s}{4}\omega$$
denote the primitive part of the Ricci form,  corresponding to the 
trace-free Ricci tensor 
$$\mathring{r}:= [r]_0= r -\frac{s}{4}g.$$ 
Suppose that 
the scalar curvature $s$ of $(M,g)$ is constant,  and set 
$$F=
\omega +  \frac{ \mathring{\rho}}{2}.$$
Then $(g,F)$ automatically solves the Einstein-Maxwell equations.
This generalizes an observation due to 
Flaherty  \citeyear{flaherty} concerning   the scalar-flat 
 case.

On a purely  calculational level, this  fact is certainly easy enough to check.
Indeed, on any oriented Riemannian $4$-manifold, 
the $2$-forms canonically decompose 
$$\Lambda^2 = \Lambda^+\oplus \Lambda^-$$
into    self-dual and anti-self-dual parts, 
and it is then easily shown that for any $2$-form
$$F=F^++F^-$$
one has 
$$\Big[ F\circ F\Big]_0= 2F^+\circ F^-.$$
Since in our special case we have 
$F^+ = \omega$ and $F^-= \mathring{\rho}/2$,
it therefore follows that 
$$\Big[ F\circ F\Big]_0= -\mathring{r},$$
so that 
$$\Big[ r+ F\circ F\Big]_0= 0$$
as required. Moreover, since $\rho$ is automatically closed,
and its self-dual part $s\omega/4$ is closed if $s$ is assumed to be constant,
we conclude that $F$  is indeed harmonic on a constant-scalar-curvature
K\"ahler surface, exactly as claimed.

But  there is actually a great deal more going on here.
Recall  that Calabi \citeyear{calabix}  defined an {\em extremal K\"ahler
metric} on a compact complex manifold $(M^{2m},J)$ to be a 
K\"ahler metric which is a critical point of the Riemannian
functional 
$$
g\longmapsto \int_Ms^2 ~d\mu
$$
restricted to a fixed K\"ahler class  $[\omega]\in H^2(M)$.
The associated Euler-Lagrange equations  then amount to requiring 
that  the gradient $\nabla s$ of the scalar curvature be the real part of a holomorphic vector field.  
In particular, any constant-scalar-curvature K\"ahler metric
is extremal in this sense. In fact, as was recently proved by 
Chen \citeyear{xxel}, extremal K\"ahler metrics actually always 
{\em minimize} $\int s^2 d\mu$ in their K\"ahler classes.
In the constant-scalar-curvature case, this was long ago shown 
by Calabi, using a simple but elegant argument. Indeed, 
if $(M^{2m}, g, J)$ is a compact K\"ahler manifold
of complex dimension $m$, then 
$$
\int_M s_g~d\mu_g= \int_M \rho \wedge \star \omega = \frac{4\pi}{(m-1)!}c_1\cdot [\omega]^{m-1}$$
and 
$$\int_M 1~d\mu = \int_M \frac{\omega^{\wedge m}}{m!} = \frac{1}{m!} [\omega]^{m}$$
so that 
the Cauchy-Schwarz inequality tells us that
\begin{equation}
\label{calabi1}
\int_Ms^2d\mu \geq \frac{(\int s~d\mu )^2}{\int 1~d\mu} = 
\frac{16\pi^2 m}{(m-1)!}\frac{(c_1\cdot [\omega]^{m-1})^2}{[\omega]^m}
\end{equation}
with equality iff $s$ is constant.

Calabi  \citeyear{calabix}
also considered 
 the Riemannian functionals 
$$
g\longmapsto \int_M |r|^2_g~d\mu_g$$
and 
$$
g\longmapsto \int_M |{\mathcal R}|^2_g~d\mu_g$$
obtained by squaring  the $L^2$-norms of the Ricci curvature $r$
and  the full Riemann curvature ${\mathcal R}$.
Here, his observation was  that the restriction of either of these
functionals to the space of K\"ahler metrics 
can be rewritten in the form 
$$g \longmapsto a + b \int_M s^2~d\mu$$
where $a$ and $b>0$ are  constants depending only
on the K\"ahler class. For example, 
$$
 \int_M |r|^2_gd\mu_g= 
 \frac{1}{2} \int_M s_g^2~d\mu_g - \frac{8\pi^2}{(m-2)!} c_1^2 \cdot [\omega]^{m-2},
$$
so that 
\begin{equation}\label{calabi2}
\int_M |r|^2_gd\mu_g \geq \frac{8\pi^2}{(m-2)!}\left[
 \frac{m}{m-1}\frac{(c_1\cdot [\omega]^{m-1})^2}{[\omega]^m}
- c_1^2 \cdot [\omega]^{m-2}\right],
\end{equation}
with equality iff $s$ is constant. 
Thus extremal K\"ahler metrics turn out to minimize these functionals, too. 
 
I would now like to point out  an interesting  Riemannian analog of 
Calabi's variational problem  that leads to the Einstein-Maxwell equations
on a smooth compact $4$-manifold. To this end, first notice that  the K\"ahler form of
a K\"ahler surface is self-dual and harmonic. Let us therefore introduce the 
following notion:

\begin{defn}
Let $M$ be  smooth compact oriented 4-manifold, and let 
$[\omega ]\in H^2(M, \RR )$ be a deRham class 
with $[\omega ]^2 > 0$. 
We will then say that a Riemannian 
metric $g$ is {\em adapted} to $[\omega ]$ if
the harmonic form $\omega$ representing
 $[\omega ]$ with respect to $g$ 
  is  {\em self-dual}.
\end{defn}

This allows us to introduce the Riemannian analog  of
a K\"ahler class: 

\begin{defn}
Let $M$ be  smooth compact oriented 4-manifold, and let 
$[\omega ]\in H^2(M, \RR )$ be a deRham class 
with $[\omega ]^2 > 0$. 
We then set
$${\mathcal G}_{[\omega ]}= \Big\{ \mbox{smooth metrics $g$ on $M$ which are 
 adapted to $[\omega ]$}\Big\}.$$
\end{defn}

In particular, if $\omega$ is the K\"ahler form of a metric $g$ on 
$M$ which is K\"ahler with respect to some complex structure $J$, then  
${\mathcal G}_{[\omega ]}$ contains the entire K\"ahler class of
$\omega$ on $(M,J)$. Of course, however, ${\mathcal G}_{[\omega ]}$ is
vastly larger than a K\"ahler class. In particular, if $g$ belongs to ${\mathcal G}_{[\omega ]}$,
so does every  conformally related metric $\tilde{g}=u^2 g$. 
It is also worth noticing  that 
$${\mathcal G}_{\lambda [\omega ]}= {\mathcal G}_{[\omega ]}$$
for any non-zero real constant $\lambda$. 

It is now  germane  to  ask precisely   how large 
${\mathcal G}_{[\omega ]}$ really 
is relative to 
$${\mathcal G}= \Big\{ \mbox{smooth metrics $g$ on $M$}\Big\}, $$
and to ponder the dependence of ${\mathcal G}_{[\omega ]}\subset {\mathcal G}$ 
on $[\omega ]$ as we allow 
this cohomology class  to vary  through the open cone
$${\mathcal C} = \Big\{ [\omega ]\in H^2 (M,\RR)~\Big|~ [\omega ]^2 > 0\Big\}.$$

\begin{prop}[Donaldson/Gay-Kirby]
Let $M$ be any smooth compact $4$-manifold with $b_+(M)\neq 0$. 
For any $[\omega ] \in {\mathcal C}$, ${\mathcal G}_{[\omega ]}\subset {\mathcal G}$
is a Fr\'echet submanifold of finite codimension $b_-(M)$. Moreover, 
${\mathcal G}_{[\omega ]}\neq \emptyset$ for all $[\omega ]$ belonging to 
an open dense subset of
${\mathcal C}$. 
\end{prop} 

Indeed, if $g\in {\mathcal G}_{[\omega ]}$ and if $\omega \in [\omega ]$ is
the harmonic representative, then Donaldson  \citeyear[p. 336]{donperiod}
has  shown that 
$T_g{\mathcal G}_{[\omega ]}$ is precisely the $L^2$-orthogonal of
the $b_-(M)$-dimensional subspace 
$$
\{ \omega \circ \varphi ~|~ \varphi \in {\mathcal H}^-_g\} \subset \Gamma (\odot^2 T^*M) , 
$$
where ${\mathcal H}^-_g$ is the space of anti-self-dual harmonic $2$-forms
with respect to $g$; moreover, his proof  also shows that the subset of 
 $[\omega]\in {\mathcal C}$ for which 
${\mathcal G}_{[\omega ]}\neq \emptyset$ is necessarily open. 
On the other hand, Gay and Kirby 
\citeyear{kirbyperiod} found an 
essentially  explicit way of constructing 
a metric $g$ adapted to 
any $[\omega ] \in {\mathcal C}\cap 
H^2 (M, {\mathbb Z})$, so that 
${\mathcal G}_{[\omega ]}\neq \emptyset$ for 
 any $[\omega ]$ in the dense  subset 
$[{\mathcal C}\cap 
H^2 (M, {\mathbb Q})]\subset {\mathcal C}$. 

We now consider the natural generalization  of Calabi's variational problem to this
broader context.

\begin{prop}
An $[\omega]$-adapted metric $g$ is a critical point 
 of the Riemannian functional 
$$
g \longmapsto \int_M s_g^2 d\mu_g
$$
restricted to ${\mathcal G}_{[\omega]}$ iff  either
\begin{itemize}
\item 
$g$ is a solution of the Einstein-Maxwell equations, in conjunction with a
unique harmonic form  $F$
with $F^+=\omega$; or else 
\item $g$ is scalar-flat $(s\equiv 0)$. 
\end{itemize}
\end{prop}

\begin{proof}
 Consider a 1-parameter family of metrics 
$$g_t:= g +th+ O(t^2)$$ 
in ${\mathcal G}_{[\omega ]}$. 
By Donaldson's result,
we know that $h$ can be taken to be 
any smooth symmetric tensor field which satisfies  
$$\int_M \langle h , \omega \circ \varphi \rangle d\mu =0$$
for all harmonic forms $\varphi \in \Gamma(\Lambda^-)$,
where $\omega$ is the $g$-harmonic representative of
$[\omega ]$.
On the other hand,  a standard calculation \cite{bes} shows that 
$$  \left.\frac{d}{dt}{s}\right|_{t=0}
 = 
\Delta h^a_a + \nabla^a\nabla^bh_{ab} -h^{ab}r_{ab},$$
and  
$$\left.\frac{d}{dt}[ d\mu ]\right|_{t=0}=\frac{1}{2}h^a_ad\mu ,$$
so that 
\bea \frac{d}{dt}\left. \left[ \int_M s^2d\mu \right]\right|_{t=0}
&=&\int_M2s\dot{s}d\mu + \int_Ms^2\dot{d\mu}\\&=&
\int 2s\left( \Delta h^a_a + \nabla^a\nabla^bh_{ab} -h^{ab}\mathring{r}_{ab} \right)d\mu
\eea
where $\mathring{r}$ again denotes the trace-free part of the Ricci tensor.

Let us now ask when a metric is critical within its conformal class.
This corresponds to setting $h=vg$ for some 
smooth function
$v$. 
We then have 
$$\frac{d}{dt}\int_M s^2d\mu = \int 2s(3\Delta v)d\mu= 
6\int_M\langle ds, dv\rangle d\mu ,$$
so the derivative is zero for all such variations
iff $s$ is constant. 

We may thus assume henceforth that $s$ is constant. We then have 
\bea \frac{d}{dt}\int_M s^2d\mu  &=&
 2s \int \left( \Delta h^a_a + \nabla^a\nabla^bh_{ab} -h^{ab}\mathring{r}_{ab} 
\right)d\mu\\&=&
-2s \int_M \langle h , \mathring{r}\rangle d\mu .\eea
If $s\equiv 0$, this obviously vanishes for every $h$, and
$g$ is a critical point. Otherwise,  $g$ will be critical 
iff $\mathring{r}$ belongs to the  $L^2$-orthogonal complement of
  $T_g{\mathcal G}_{[\omega ]}$.
But we already have seen
 that this orthogonal complement precisely consists of  tensors of the form $\omega \circ \varphi$,
$\varphi \in {\mathcal H}^-_g$. Thus, when $s\not\equiv 0$, 
$g$ is a critical point iff $s$ is constant and 
$\mathring{r} = \omega \circ \varphi$  for some $\varphi \in {\mathcal H}^-_g$. 
But, setting
$$F= \omega + \frac{\varphi}{2}~,$$
this  is in turn  equivalent to saying that $(g, F)$ 
solves  the Einstein-Maxwell equations, as claimed.  
\end{proof}
%

%


So, why are constant-scalar-curvature K\"ahler metrics critical points 
of $\int s^2d\mu$ restricted to ${\mathcal G}_{[\omega ]}$?
Well, we will now see that  they typically turn out not only to be critical points, 
but  actually to be {\em minima}. Indeed, the following
result \cite{lpm,lric,lebsurv2} 
may be thought of as a Riemannian generalization  of 
 Calabi's inequalities (\ref{calabi1}--\ref{calabi2}): 

\begin{thm}
\label{whatsit} 
Let $(M^4,J)$ be a compact complex surface, and suppose that 
$[\omega]$ is a K\"ahler class with $c_1\cdot [\omega] \leq 0$. Then 
any {\em Riemannian} metric $g\in {\mathcal G}_{[\omega ]}$ satisfies
the inequalities
\begin{eqnarray}
\int s^2 d\mu &\geq& 32\pi^2 \frac{(c_1\cdot [\omega ])^2}{[\omega ]^2}
\label{who1}
\\
\int |r|^2 d\mu &\geq& 8\pi^2 \Big[ 2\frac{(c_1\cdot [\omega ])^2}{[\omega ]^2} -c_1^2\Big]
\label{who2}
\end{eqnarray}
with equality if and only if $g$ is  constant-scalar-curvature 
K\"ahler.
\end{thm} 

In the equality case, the complex structure $\tilde{J}$ with respect to which 
$g$ is K\"ahler will typically be different from $J$, but must have the same
first Chern class $c_1$, while  its K\"ahler class must be a positive multiple of 
$[\omega ]$.

We also remark that if $(M,J)$ is not rational or ruled, the hypothesis 
that $c_1\cdot [\omega ] \leq 0$ holds automatically, and that
in this setting 
 a K\"ahler metric is extremal iff it has constant
scalar curvature. In this context, the relevant 
constant is of course necessarily non-positive. 

By contrast, if $(M,J)$ is rational or ruled, there will
always be K\"ahler classes for which $c_1\cdot [\omega ] > 0$.
When this happens, 
the above generalization (\ref{who1}) of (\ref{calabi1}) 
  turns out   definitely {\em not}  to 
hold for arbitrary Riemannian metrics. 
Instead, the correct generalization \cite{lcp2}  is that 
\begin{equation}\label{who3}
Y_{[g]} \leq    \frac{4\pi ~c_1\cdot [\omega ]}{\sqrt{[\omega ]^2/2}}~,
\end{equation}
where the Yamabe constant $Y_{[g]}$ is obtained by minimizing the 
Einstein-Hilbert action $\int s~d\mu$ over all unit-volume
metrics $\tilde{g}=u^2g$ conformal to $g$. Moreover, the inequality 
is strict unless the Yamabe minimizer is a constant-scalar-curvature 
K\"ahler metric, so  that (\ref{who1}) 
is in fact violated by an appropriate  conformal rescaling of any generic  
 Riemannian metric of positive scalar curvature. 

It is also worth remarking that no sharp lower bound in the spirit of 
 Theorem \ref{whatsit} is currently known for the square-norm 
 $\int |{\mathcal R}|^2d\mu$ of the Riemann curvature tensor. 
Deriving  one would be extremely interesting and potentially very useful,
but, for reasons  I will now explain, the technical obstacles
to doing so  seem formidable.

Recall that, by raising an index, the 
Riemann curvature tensor may be reinterpreted as a linear map 
$\Lambda^2 \to \Lambda^2$,  called the {\em curvature operator}.
The 
decomposition $\Lambda^2 = \Lambda^+\oplus \Lambda^-$ 
thus  allows one to view  this linear map  as consisting of four blocks: 
\begin{equation}
\label{curv}
{\mathcal R}=
\left(
\mbox{
\begin{tabular}{c|c}
&\\
$W_++\frac{s}{12}$&$\mathring{r}$\\ &\\
\cline{1-2}&\\
$\mathring{r}$ & $W_-+\frac{s}{12}$\\&\\
\end{tabular}
} \right) . 
\end{equation}
Here $W_\pm$ are the trace-free pieces of the appropriate blocks,
and  are  called the
self-dual and anti-self-dual Weyl curvatures, respectively. 
The scalar curvature  $s$ is understood to act by scalar multiplication,
whereas the   
     trace-free Ricci curvature
$\mathring{r}=r-\frac{s}{4}g$ 
acts on 2-forms by
$\varphi_{ab} \mapsto ~ 
2{\varphi_{[a}}^{c} \mathring{r}_{b]c}$.

When $(M,g)$ happens to be K\"ahler, ${\Lambda^{2,0}} \subset \ker {\mathcal R}$,
and the entire upper-left-hand block is therefore entirely determined by 
the scalar curvature $s$. For K\"ahler metrics, one thus obtains the identity 
$$|W_+|^2 \equiv \frac{s^2}{24},$$
and  Gauss-Bonnet-type formul{\ae}
 like
$$
 (2\chi + 3\tau ) (M)  
  = \frac{1}{4\pi^2} \int_M \left(\frac{s^2}{24}+ 2 |W_+|^2 - \frac{|\mathring{r}|^2}{2}\right)d\mu  
$$
reduce many questions about square-norms of curvature to questions 
about the scalar-curvature alone. But 
for general Riemannian metrics, the norms of $s$ and $W_+$ are
utterly independent quantities, so if one wants  to use the identity 
\begin{equation}\label{whyzit} 
\int|r|^2d\mu = -8\pi^2 (2\chi + 3\tau ) (M)+ 8 \int\left(\frac{s^2}{24}+ \frac{1}{2}
 |W_+|^2 \right)d\mu  
\end{equation}
to prove a generalization of  (\ref{calabi2}) for Riemannian metrics, 
information must be obtained concerning  not only the scalar curvature, 
but also concerning the self-dual Weyl curvature as well.

The curvature estimates  of Theorem \ref{whatsit} 
are derived by means of  Seiberg-Witten theory \cite{witten}, 
making it clear that this really is an essentially
$4$-dimensional story. The complex structure $J$ determines
a spin$^c$ structure on $M$ with 
twisted spin bundles ${\mathbb S}_\pm \otimes L^{1/2}$,
where $L^{-1}$ is the  canonical line bundle $\Lambda^{2,0}$
of $(M,J)$.  For simplicity, suppose  that $c_1\cdot [\omega ]< 0$. 
For each metric $g\in {\mathcal G}_{[\omega ]}$,
one  then considers the Seiberg-Witten equations 
\begin{eqnarray*} D_{A}{\Phi}&=&0\\
 {F}^+_{A}&=&-\frac{1}{2}{\Phi}\odot {\bar{\Phi} } \end{eqnarray*}
 where the unknowns are 
 a unitary connection $A$ on the line-bundle  $L\to M$
 and a twisted spinor
 $\Phi \in \Gamma ({\mathbb S}_+\otimes L^{1/2})$;
 here 
 $D_A:  \Gamma ({\mathbb S}_+\otimes L^{1/2})\to  
 \Gamma ({\mathbb S}_-\otimes L^{1/2}) $ denotes the twisted Dirac operator 
 associated with $A$, and $F^+_A$ is the self-dual part of the 
 curvature of $A$. 
 One then shows that there must be at least one solution for each $g\in {\mathcal G}_{[\omega ]}$ 
 by establishing 
 a count of solutions modulo gauge equivalence which is independent
 of the metric and which is obviously non-zero for  a K\"ahler metric.

However, the Seiberg-Witten equations can be shown 
to imply various curvature estimates via Weitzenb\"ock formul{\ae}. 
In particular, the existence of at least one solution for each
metric $\tilde{g}=u^2 g$ conformal to $g$ is enough to guarantee that 
the curvature of $g$ satisfies
\begin{eqnarray*}
\int_{M} {s}^2
{d\mu}_{g}& \geq &32\pi^2[c_1 (L)^+]^2
\\
\int_{M} \left({s}-\sqrt{6}|{W}_+|\right)^2 {d\mu}_{g}
 &\geq &  72\pi^2[ c_1(L)^+]^2
\end{eqnarray*}
where $[c_1(L)]^+$ is the orthogonal projection of 
$c_1(L)\in H^2(M, \RR) = {\mathcal H}^+_g\oplus  {\mathcal H}^-_g$
into the space ${\mathcal H}^+_g$ of harmonic self-dual $2$-forms, 
defined with respect to $g$. 
Since $\omega$ is assumed to be self-dual with respect to $g$,
we therefore have 
$$
[c_1 (L)^+]^2 \geq \frac{(c_1\cdot [\omega ])^2}{[\omega ]^2}
$$
by the Cauchy-Schwarz inequality. 
Inequality (\ref{who1}) follows. 
Since yet another Cauchy-Schwarz argument shows that 
$$
 \int\left(\frac{s^2}{24}+ \frac{1}{2}
 |W_+|^2 \right)d\mu_g  
\geq \frac{1}{36} \int \left({s}-\sqrt{6}|{W}_+|\right)^2 {d\mu}_{g}
$$
the second inequality and (\ref{whyzit}) together imply
(\ref{who2}).  The fact that only K\"ahler metrics can saturate 
(\ref{who1}) or (\ref{who2}) is then deduced by examining the 
relevant Weitzenb\"ock formul{\ae}.

One might be tempted to expect the story to be similar 
for the norm of the full Riemann tensor. 
After all, the identity 
$$
\int_M |{\mathcal R}|^2 d\mu = -8\pi^2 (\chi + 3\tau )(M) +2\int_M \left(\frac{s^2}{24}+ 2
 |W_+|^2 \right)d\mu_g  
$$
 certainly provides a good analog  of (\ref{whyzit}). 
In the K\"ahler case, one has 
$$\frac{s^2}{24} = |W_+|^2,$$
so  this simplifies to become
$$
\int_M |{\mathcal R}|^2 d\mu = 8\pi^2 (c_2-c_1^2) +\frac{1}{4}\int_M s^2 ~d\mu_g  ~,
$$
and  applying (\ref{calabi1}) we therefore obtain  Calabi's inequality 
\begin{equation}
\label{calabi4}
\int_M |{\mathcal R}|^2 d\mu\geq 8\pi^2  \left[ \frac{(c_1\cdot [\omega ])^2}{[\omega ]^2} 
+ c_2 - c_1^2  \right] 
\end{equation}
for any K\"ahler metric. In light of Theorem \ref{whatsit}, 
it might therefore seem reasonable to hope that one could simply 
extend this inequality to general Riemannian metrics by means
of Seiberg-Witten theory. However, we
will now show that  this cannot work. 
The key idea is to examine certain extremal K\"ahler metrics 
from the vantage point of their {\em reversed orientations}.

\begin{thm}
Calabi's inequality (\ref{calabi4}) cannot possibly be extended to general 
Riemannian metrics by means of Seiberg-Witten theory. Indeed, 
there actually exist smooth compact 
oriented Riemannian $4$-manifolds $(M,g)$ which admit a spin$^c$ structure of
almost-complex type with non-zero Seiberg-Witten invariant, but  such that 
$$
\int_M |{\mathcal R}|^2 d\mu < 8\pi^2  \left[ \frac{(c_1\cdot [\omega ])^2}{[\omega ]^2} 
-(\chi + 3\tau )(M)  \right] 
$$
for some self-dual harmonic $2$-form $\omega$ on $(M,g)$. 
\end{thm}
  \begin{proof}
 A  {\em Kodaira-fibered complex surface} is by definition 
 a  compact complex surface $X$ equipped with 
   a holomorphic
submersion $\varpi: X\to \mathscr{B}$  onto a  
compact complex curve, such that the  base $\mathscr{B}$ and fiber 
$\mathscr{F}_z= \varpi^{-1}(z)$ both have genus $\geq 2$. 
The product $\mathscr{B}\times \mathscr{F}$ of two complex curves of
genus $\geq 2$ is certainly  Kodaira fibered, but such  
a product   signature
$\tau=0$. However,
one can  also construct examples \cite{atkodf,kodf}
with $\tau > 0$ by taking {\em branched covers} of 
products. 

Let $X$ be  any such Kodaira-fibered surface
with $\tau (X)> 0$, and let  
$\varpi : X\to \mathscr{B}$ be its Kodaira fibration.
Let $p$ denote the the genus of the base $\mathscr{B}$, and let
$q$ denote the genus of some fiber $\mathscr{F}$ of $\varpi$.
A beautiful result of Fine \citeyear{fine}  then 
asserts that $X$  actually  
admits a family of extremal  K\"ahler metrics;
namely, {\em for any sufficiently small $\epsilon > 0$, 
$$
[\omega_{\epsilon}] = 2(p-1) \mathscr{F}- \epsilon c_{1}
$$
is a K\"ahler class on $X$ which is represented by 
a K\"ahler metric $g_\epsilon$ of constant scalar curvature.}

These metrics, being K\"ahler,  have total scalar curvature 
$$\int s_{g_{\epsilon}} d\mu_{g_{\epsilon}} = 4\pi c_{1}\cdot 
[\omega_{\epsilon}]= -4\pi (\chi +\epsilon c_{1}^{2}) (X)$$
and total volume
$$\int d\mu_{g_{\epsilon}} = \frac{[\omega_{\epsilon}]^{2}}{2}=
\frac{\epsilon}{2}(2\chi +\epsilon c_{1}^{2}) (X).$$
Since $s_{g_{\epsilon}}$ is constant, 
it follows  that
\begin{eqnarray*}
\int s^{2}_{g_{\epsilon}}d\mu_{g_{\epsilon}}	 & = & 
\frac{32\pi^{2}}{\epsilon} 
\frac{(\chi +\epsilon c_{1}^{2})^{2}}{2\chi +\epsilon c_{1}^{2}} . 
\end{eqnarray*}
These metrics therefore satisfy 
\begin{eqnarray*}
\int_X |{\mathcal R}|^2_{g_\epsilon} d\mu_{g_\epsilon}
 &=& 8\pi^2 (c_2-c_1^2) +\frac{1}{4}\int_X s^2 ~d\mu_g \\
&=& 8\pi^2 \left[ -(\chi + 3\tau)(X) +
\frac{(\chi +\epsilon c_{1}^{2})^{2}}{\epsilon(2\chi +\epsilon c_{1}^{2})}\right] 
\end{eqnarray*}

On the other hand, there are symplectic forms on $X$ which 
are compatible with the {\em non-standard} orientation 
of $X$; for example, the cohomology class $\mathscr{F}+\varepsilon c_{1}$
is represented by such forms if $\varepsilon$ is sufficiently small. 
A celebrated theorem of Taubes \citeyear{taubes} therefore 
tells us that the reverse-oriented version
$M=\overline{X}$ of $X$ has a non-trivial Seiberg-Witten 
invariant  \cite{leung,kot2}. The relevant 
spin$^{c}$ structure on $\overline{X}$ is of 
almost-complex type, and its first Chern class, which
we will denote by $\bar{c}_{1}$, is given by
$$\bar{c}_{1}= c_{1}+4(p-1)\mathscr{F}.$$
Since this a $(1,1)$-class, 
  one has
$$
({\bar{c}_1})^{+}	=
\frac{{\bar{c}_1}\cdot [\omega_{\epsilon}]}{[\omega_{\epsilon}]^{2}}
\omega_{\epsilon}  
	  =  -\frac{(\chi + 3\epsilon \tau)}{[\omega_{\epsilon}]^{2}}
	 \omega_{\epsilon} ,
$$
relative to the K\"ahler metric 
$g_{\epsilon}$, so that 
$$
|({\bar{c}_1})^{+}|^{2}	  =  
\frac{(\chi + 3\epsilon \tau)^{2}}{[\omega_{\epsilon}]^{2}} =
\frac{(\chi + 3\epsilon \tau)^{2}}{\epsilon(2\chi +\epsilon c_{1}^{2})} .
$$
Now since ${\bar{c}_1}$ arises from an almost-complex structure on 
$\overline{X}$, we have 
$$|({\bar{c}_1})^{-}|^{2}-|({\bar{c}_1})^{+}|^{2}= 2\chi - 3\tau, $$
so that 
$$
|({\bar{c}_1})^{-}|^{2}	  =  2\chi - 3\tau  +
 \frac{(\chi + 3\epsilon \tau)^{2}}{\epsilon(2\chi +\epsilon c_{1}^{2})}, 
$$
and 
$$
|({\bar{c}_1})^{-}|^{2}	- (\chi - 3\tau )(X) = \chi (X) +  \frac{(\chi + 3\epsilon \tau)^{2}}{\epsilon(2\chi +\epsilon c_{1}^{2})},
$$
where, for example,    $\tau$  indicates $\tau (X)$. But it therefore  follows  that 
\begin{eqnarray*}
\frac{1}{8\pi^2}\int_X |{\mathcal R}|^2_{g_\epsilon} d\mu_{g_\epsilon}
- \Big[ |({\bar{c}_1})^{-}|^{2}	- (\chi - 3\tau )\Big]
 &=&  -(2\chi + 3\tau) + 
\frac{(\chi +\epsilon c_{1}^{2})^{2}- (\chi + 3\epsilon \tau)^{2}}{\epsilon(2\chi +\epsilon c_{1}^{2})}
 \\&=& -(2\chi + 3\tau) + 
\frac{4\chi \epsilon(\chi + 3\epsilon \tau) + 4\chi^2\epsilon^2}{\epsilon(2\chi +\epsilon c_{1}^{2})}
\\&=& -(2\chi + 3\tau) + 2\chi
\frac{2\chi  +\epsilon c_{1}^{2}+ 3\epsilon \tau}{2\chi +\epsilon c_{1}^{2}}
 \\&=& -3 \tau (X) \left[ 1-   \frac{2\chi  \epsilon}{2\chi + \epsilon c_1^2} \right] ~,
\end{eqnarray*}
which  is negative for any sufficiently small $\epsilon$. 
The result therefore follows once we 
take ``$\omega$'' to be 
the anti-self-dual harmonic form $({\bar{c}_1})^-$, which becomes
self-dual on $M=\overline{X}$. 
\end{proof}

Similarly,  careful examination of 
these examples also shows that, for any constant $t> 1$,   the Seiberg-Witten 
equations cannot imply an estimate of 
$$
\int \left( s- t \sqrt{6} |W_+|\right)^2 d\mu~
$$
which is saturated by constant-scalar-curvature K\"ahler metrics. 
Of course, the Seiberg-Witten equations still imply lower 
bounds for such quantities, but they  are simply never as sharp as those 
obtained for  $t\in [0,1]$.

In this article, we have  seen that constant-scalar-curvature K\"ahler metrics 
 occupy a privileged position in $4$-dimensional Riemannian
 geometry. I would therefore like  to conclude this discussion by indicating
  a bit of what we now know
 concerning  their existence. 

There are several ways to phrase the problem. From the Riemannian
point of view, one might want to fix a smooth compact oriented
$4$-manifold $M$, and simply ask whether there exists
an extremal K\"ahler metric $g$, where the
associated complex structure is not specified as part of the problem. 
Since $M$ must in particular admit a K\"ahler metric, two necessary conditions 
are  that $M$ must admit
a complex structure and have  even first Betti number.
Provided these desiderata are fulfilled, 
\citeasnoun{yujen} has then shown
that  an extremal  metric $g$  always exists.
For all but two diffeotypes, moreover, one can actually  arrange for the 
extremal K\"ahler metric $g$ 
 to have constant scalar curvature. However, these two exceptional diffeotypes
 are $\CP_2\#\overline{\CP}_2$, and $\CP_2\#2\overline{\CP}_2$, and it is now 
 known \cite{chenlebweb} that both these manifolds carry Einstein metrics ---
 indeed, even 
 Einstein metrics which are conformal rescalings of  extremal K\"ahler metrics! 
 Since, in conjunction with   $F=0$,  
 any  Einstein  metric of course satisfies the Einstein-Maxwell equations,
 we thus immediately deduce  the following: 
 
 \begin{thm} \label{freebie} 
 Let $M$ be the underlying $4$-manifold of any 
 compact complex surface of K\"ahler type. Then $M$ admits
 a Riemannian solution $(g,F)$ of
 the Einstein-Maxwell equations.  
 \end{thm}
 
 While the above formulation of Shu's result   certainly suffices 
 to   imply  Theorem \ref{freebie},
  it unfortunately also    obfuscates the 
nature  of the proof, which involves
 constructing   extremal K\"ahler metrics 
compatible with some {\em fixed} complex structure in  each possible 
deformation class. The key 
tool used for this purpose   is  due to 
Arezzo and Pacard  
\citeyear{arpa1,arpa2}, 
 who have shown that constant-scalar-curvature
K\"ahler metrics can be constructed on blow-ups and desingularizations
 of constant-scalar-curvature
K\"ahler orbifolds, under only very mild assumptions on the complex automorphism
group; similar results moreover have even been proved concerning 
  the strictly extremal case \cite{arpasing}. These gluing results 
represents a vast  generalization of earlier work
 by  the present author and his collaborators \cite{klp,leblown,lebsing1}  
regarding the limited realm 
of scalar-flat K\"ahler surfaces. 
In fact, by invoking the theory of  K\"ahler-Einstein metrics \cite{aubin,yau}, 
\citeasnoun{arpa1} had already shown 
 that  every K\"ahler-type  complex surface  of Kodaira dimension $0$ or $2$
admits compatible
constant-scalar-curvature K\"ahler metrics.  Shu's 
results concerning the  remaining cases of
  Kodaira dimensions $-\infty$ and $1$  
  are  much less robust, 
  but  still easily produce enough examples to  imply most of Theorem \ref{freebie}. 
    
 Of course, one  ultimately doesn't  want to settle for  mere existence statements;
 we  would really like to completely understand the moduli space of solutions! 
 From this 
 point of view, the first  question  to ask is whether
  there
 can only be one solution for  any given  complex structure and  K\"ahler class. Modulo
 complex automorphisms, 
uniqueness always holds in this setting, 
 as was proved in a series of a 
 fundamental papers by Donaldson 
  \citeyear{donaldsonk1}, Mabuchi \citeyear{mabuni},  and Chen \& Tian \citeyear{xxgang}. 
 For a fixed complex structure,  one also knows   that  the K\"ahler classes of  
 extremal K\"ahler metrics    sweep out an open subset of the K\"ahler cone
 \cite{ls2}, and 
  somewhat weaker results are also available regarding deformations of 
 complex structure \cite{fusch,ls}. However, it turns out that 
 the set of K\"ahler classes which are representable by extremal
 K\"ahler metrics may sometimes  be a {\em proper} non-empty  open subset 
 of the K\"ahler cone \cite{apoto,ross}. The latter phenomenon is 
 related to algebro-geometric stability problems \cite{mabdon,rotho} 
 in a manner
 which is still only partly understood, but  there  is  reason to hope   that 
 a definitive understanding of such issues may  result 
 from the incredible ferment of research  currently being carried out 
  the field. 
 
As we saw in Theorem \ref{freebie}, 
 K\"ahler geometry supplies a 
 natural and beautiful way of constructing solutions of
 the Einstein-Maxwell equations on many compact 
 $4$-manifolds. In the opposite direction,  we also have the 
 following easy but rather suggestive result:

 \begin{prop} \label{curio} 
 Let $M$ be the underlying $4$-manifold of any 
 compact complex surface of {\em non-K\"ahler}
 type with vanishing geometric genus. Then $M$
does not carry   any Riemannian 
 solution of the Einstein-Maxwell 
 equations. 
 \end{prop} 
 \begin{proof}
Let us begin by remembering   the  remarkable fact
 \cite{siu,bpv,nick} that a compact complex surface is of K\"ahler 
 type
 iff it has $b_1$ even. Consequently, for  any 
  non-K\"ahler-type complex surface $M$, $b_1$ is odd, and 
 $b_+=2p_{\bf g}$, where  $p_{\bf g}=h^{2,0}$ is  the geometric
 genus \cite[Theorem IV.2.6]{bpv}.
  Since the latter is assumed  to 
 vanish, $M$ then has negative-definite intersection form, and 
 Hodge theory tells us that there are no non-trivial 
 self-dual  harmonic $2$-forms for any metric $g$ on $M$. 
 
 Now suppose that  $(g,F)$ is  a Riemannian solution of the 
 Einstein-Maxwell equations on $M$. 
 Then  the harmonic $2$-form $F$ satisfies  $F^+=0$,
 and hence 
 $$\mathring{r}= -[F\circ F]_0=-2F^+\circ F^-=0.$$ 
 The metric 
  $g$  must therefore  be Einstein.  But the negative-definiteness 
   of the intersection form also tells us that $(2\chi + 3\tau )(M) = c_1^2 (M) \leq 0$. 
The Hitchin-Thorpe inequality for Einstein manifolds
\cite{hit} therefore guarantees that $M$ has a finite normal cover $\tilde{M}$ 
which is 
diffeomorphic to either  $K3$ or $T^4$,  and so, in particular, has 
$b_1$  even. 
Pulling back the complex structure $J$ of $M$ to this cover, 
we therefore obtain a complex surface $(\tilde{M}, \tilde{J})$ 
of K\"ahler type. 
 Averaging an arbitrary K\"ahler metric $h$ on $(\tilde{M}, \tilde{J})$  
over the finite group of deck transformation of $\tilde{M}\to M$ then gives us a K\"ahler metric
which descends to $(M,J)$. Thus $(M,J)$ is of K\"ahler type, in contradiction 
to our hypotheses. Our supposition was therefore false, and 
 $M$ thus cannot  carry any  Riemannian solution of the
Einstein-Maxwell equations. 
 \end{proof}

This article has  endeavored to  convince the reader that 
 the four-dimensional Einstein-Maxwell equations 
 represent a beautiful and natural  generalization of
 the constant-scalar-curvature K\"ahler condition. 
However, it still remains to be seen whether 
solutions   exist  on many compact $4$-manifolds  
other than  complex surfaces of K\"ahler type. 
 In this direction, my 
guess  is that Proposition \ref{curio} will  actually   prove to be 
 rather misleading. 
 For example,  there are ALE  Riemannian solutions of the 
  Einstein-Maxwell equations,  constructed \cite{yuille} via  the 
   Israel-Wilson ansatz,    
 on manifolds  very  unlike any complex surface. 
   It thus  seems reasonable to conjecture
    that there are plenty  of  
of compact solutions that in no sense  arise from  K\"ahler geometry. 
Perhaps    some interested reader will feel inspired to 
 go out and find some!  
 
 \bigskip
 
 \bigskip

 \noindent
 {\bf Acknowledgments.} The author would like to thank 
 Maciej Dunajski,  Caner Koca,  and Galliano Valent 
 for their interesting comments and queries regarding an earlier  version of this article,
as  some of these  resulted in  important  improvements to the paper.

  \end{document}